\input amstex
\input amsppt.sty
\NoBlackBoxes
\def\Im{\operatorname{Im}}
\def\Re{\operatorname{Re}}
\def\rank{\operatorname{rank}}
\def\sign{\operatorname{sign}}
\def\om{\omega}
\def\C{\Bbb C}
\def\del{\partial}
\def\w{\omega}
\def\s{\sigma}
\def\R{\Bbb R}
\def\N{\Bbb N}

\topmatter
\title
$3$-forms and almost complex structures on $6$-dimensional manifolds
\endtitle
\author
Martin Pan\'ak
\footnote{\quad Supported by the Grant
Agency of  the Czech Republic,
grant No. 201/02/1390.\hfill\hfill}, Ji\v r\'\i\ Van\v zura
\footnote{\quad Supported by the Grant Agency of Czech Academy
of Sciences, grant No. A 101 9204.\hfill\hfill}
\endauthor
\endsubjclass
\rightheadtext{$3$-forms on $6$-dimensional manifolds}
\leftheadtext{Martin Pan\'ak, Ji\v r\'\i\ Van\v zura}
\abstract
This article deals with $3$-forms on $6$-dimensional manifolds, the first
dimension where the classification of $3$-forms is not trivial.
There are three classes of multisymplectic $3$-forms there. 
We study the class which is closely related to almost
complex structures.
\endabstract
\endtopmatter

Let $V$ be a~real vector space. Recall that a~$k$-form $\omega$
($k\geq2$) is called multisymplectic if the homomorphism
$$
\iota:V\rightarrow\Lambda^{k-1}V^*,\quad v\mapsto\iota_v\omega=\omega(v,\dots)
$$
is injective.
There is a~natural action of
the general linear group $G(V)$
on $\Lambda^kV^*$, and also on  $\Lambda_{ms}^kV^*$,
the subset of the multisymplectic forms. Two multisymplectic
forms are called equivalent if they belong to the same orbit of the action.
For any form
$\omega\in\Lambda^kV^*$ define a~subset
$$
\Delta(\omega)=\{v\in V;(\iota_v\omega)\wedge(\iota_v\omega)=0\}.
$$

If $\dim V=6$ and $k=3$ the subset $\Lambda_{ms}^3V^*$ consists of three
orbits. Let $e_1,\dots,e_6$ be a basis of $V$ and
$\alpha_1,\dots,\alpha_6$ the corresponding dual basis. Representatives
of the three orbits can be expressed in the form
\roster
\item"{(1)}"$\omega_1=\alpha_1\wedge\alpha_2\wedge\alpha_3+
           \alpha_4\wedge\alpha_5\wedge\alpha_6$,
\item"{(2)}"$\omega_2=\alpha_1\wedge\alpha_2\wedge\alpha_3+
             \alpha_1\wedge\alpha_4\wedge\alpha_5+
             \alpha_2\wedge\alpha_4\wedge\alpha_6-
             \alpha_3\wedge\alpha_5\wedge\alpha_6$,
\item"{(3)}"$\omega_3=\alpha_1\wedge\alpha_4\wedge\alpha_5+
             \alpha_2\wedge\alpha_4\wedge\alpha_6+
             \alpha_3\wedge\alpha_5\wedge\alpha_6$.
\endroster
Multisymplectic forms are called of type 1, resp. of type 2, resp. of type 3
accordning to which orbit they belong to.
There is the following characterisation of the orbits:
\roster
\item"{(1)}" $\om$ is of type 1 if and only if $\Delta(\omega)=V^a\cup V^b$,
where $V^a$ and $V^b$ are 3-dimensional subspaces satisfying $V^a\cap
V^b=\{0\}$.
\item"{(2)}" $\om$ is of type 2 if and only if $\Delta(\omega)=\{0\}$.
\item"{(3)}" $\om$ is of type 3 if and only if $\Delta(\omega)$ is a
3-dimensional subspace.
\endroster

The forms $\omega_1$ and $\omega_2$ have equivalent
complexifications. From this point of view the forms of type 3
are exceptional. You can find more  about these forms
in [V].

A multisymplectic $k$-form on a~manifold $M$
is a~section of $\Lambda^kT^*M$ such that its restriction to the tangent
space $T_xM$ is multisymplectic for any $x\in M$,
and is of type $i$ in $x\in M$, $i=1$, $2$,
$3$, if the restriction to $T_xM$ is of type $i$.
A multisymplectic form
on $M$ can change its type as seen on
$$
\split
\sigma=dx_1\wedge dx_2\wedge dx_3 + dx_1\wedge dx_4\wedge dx_5 + dx_2\wedge dx_4\wedge dx_6 +
\cr\sin(x_3 + x_4) dx_3\wedge dx_5\wedge dx_6 +
\sin(x_3 + x_4) dx_4\wedge dx_5\wedge dx_6,
\endsplit
$$
a $3$-form on $\R^6$.
$\sigma$ is of type $3$ on the submanifold
given by the equation $x_3+x_4=k\pi$, $k\in \N$.
If $x_3+x_4\in(k\pi,(k+1)\pi)$, $k$ even, then
$\s$ is of type $1$ and if $x_3+x_4\in(k\pi,(k+1)\pi)$, $k$ odd,
then $\s$ is of type $2$. Let us point out that $\s$ is closed
and  invariant under the action of the group
$(2\pi\Bbb Z)^6$ and we can factor $\s$ to get a~form changing the
type on $\Bbb R^6/(2\pi\Bbb Z)^6$, which is the $6$-dimensional
torus, that is  $\s$ is closed on a~compact manifold.

The goal of this paper is to study the forms of type 2. We
denote $\omega=\omega_2$.
\head
\centerline{\bf $3$-forms on vector spaces}
\endhead
Let $J$ be an automorphism of  a~6-dimensional real vector space $V$
satisfying $J^2=-I$. Further let $V^\C=V\oplus iV$ be the
complexification of $V$. There is the standard decomposition $V^\C=
V^{1,0}\oplus V^{0,1}$. Consider  a~non-zero form $\gamma$
of type $(3,0)$ on $V^\C$ and  set
$$
\gamma_0=\Re\gamma,\quad\gamma_1=\Im\gamma.
$$
For any $v_1\in V$ there is $v_1+iJv_1\in V^{0,1}$, and consequently
 $\gamma(i(v_1+iJv_1),v_2,v_3)=0$ for any $v_2,v_3\in V$. This implies
$\gamma_0(i(v_1+iJv_1),v_2,v_3)=0$ and  $\gamma_1(i(v_1+iJv_1),v_2,v_3)=0$.
Thus
$$
0=\gamma_0(i(v_1+iJv_1),v_2,v_3)=\gamma_0(iv_1,v_2,v_3)-\gamma_0(Jv_1,v_2,v_3).
$$

\noindent
Similarly we can proceed with $\gamma_1$ and we get
$$
\gamma_0(iv_1,v_2,v_3)=\gamma_0(Jv_1,v_2,v_3),\quad
\gamma_1(iv_1,v_2,v_3)=\gamma_1(Jv_1,v_2,v_3)
$$
for any $v_1,v_2,v_3\in V$. Moreover there is
$$
\split
\gamma_0(w_1,w_2,w_3)&=\Re(-\gamma(i^2w_1,w_2,w_3))=\Re(-i\gamma(iw_1,w_2,w_3))\\
&=\Im(\gamma(iw_1,w_2,w_3))=\gamma_1(iw_1,w_2,w_3),
\endsplit
$$
for any $w_1$, $w_2$, $w_3\in V^\C$ and that is $\gamma_1(w_1,w_2,w_3)=
-\gamma_0(iw_1,w_2,w_3)$.
Finally,
$$
\split
\gamma_0(Jv_1,v_2,v_3)&=\gamma_0(iv_1,v_2,v_3)=\Re(\gamma(iv_1,v_2,v_3))=
\Re(i\gamma(v_1,v_2,v_3))\\
&=\Re(\gamma(v_1,iv_2,v_3))=\Re(\gamma(v_1,Jv_2,v_3))=\gamma_0(v_1,Jv_2,v_3).\\
\endsplit
$$

\noindent
Along these lines we obtain
$$
\split
\gamma_0(Jv_1,v_2,v_3)&=\gamma_0(v_1,Jv_2,v_3)=\gamma_0(v_1,v_2,Jv_3),\\
\gamma_1(Jv_1,v_2,v_3)&=\gamma_1(v_1,Jv_2,v_3)=\gamma_1(v_1,v_2,Jv_3),
\endsplit
$$
that is both forms $\gamma_0$ and $\gamma_1$ are pure with
respect to the complex structure $J$.

\proclaim{1. Lemma}
The real 3-forms $\gamma_0|V$ and $\gamma_1|V$ (on $V$) are multisymplectic.
\endproclaim
\demo{Proof}
Let us assume that $v_1\in V$ is a~vector such that for any vectors
$v_2,v_3\in V$  $(\gamma_0|V)(v_1,v_2,v_3)=0$ or equivalently
$\gamma_0(v_1,v_2,v_3)=0$.
There are uniquely determined vectors $w_1,w_2,w_3\in V^{1,0}$ such that
$$
v_1=w_1+\bar{w}_1,\quad v_2=w_2+\bar{w}_2,\quad v_3=w_3+\bar{w}_3.
$$
Then
$$
\split
0=\gamma_0(v_1,v_2,v_3)&=\Re(\gamma(w_1+\bar{w}_1,w_2+\bar{w}_2,
w_3+\bar{w}_3))\\
&=\Re(\gamma(w_1,w_2,w_3))=\gamma_0(w_1,w_2,w_3)
\endsplit
$$
(for a~fixed $w_1$, and arbitrary $w_2,w_3\in V^{1,0}$).
Because $iw_2\in V^{1,0}$, we find that
$$
\gamma_0(iw_1,w_2,w_3)=\gamma_0(w_1,iw_2,w_3)=0.
$$
Moreover $\gamma_1(w,w',w'')=-\gamma_0(iw,w',w'')$ for any
$w,w',w''\in V^\C$, and we get
$$
\gamma_1(w_1,w_2,w_3)=-\gamma_0(iw_1,w_2,w_3)=0
$$
for arbitrary $w_2,w_3\in V^{1,0}$. Thus
$$
\gamma(w_1,w_2,w_3)=\gamma_0(w_1,w_2,w_3)+i\gamma_1(w_1,w_2,w_3)=0
$$
for arbitrary $w_2,w_3\in V^{1,0}$.

Because $\gamma$ is a~non-zero
complex 3-form on the complex 3-dimensional vector space $V^{1,0}$, we
find that $w_1=0$, and consequently $v_1=0$. This proves that the real
3-form $\gamma_0|V$ is multisymplectic.
We find that the real 3-form
$\gamma_1|V$ is also multisymplectic likewise.
\enddemo

\proclaim{2. Lemma}
The forms $\gamma_0|V$ and $\gamma_1|V$ satisfy $\Delta(\gamma_0|V)=
\{0\}$ and $\Delta(\gamma_1|V)=\{0\}$.
\endproclaim
\demo{Proof}
The complex 3-form $\gamma$ is decomposable, and therefore
$\gamma\wedge\gamma=0$. This implies that for any $w\in V^\C$
$(\iota_w\gamma)\wedge(\iota_w\gamma)=0$. Similarly  for any $w\in
V^\C$ $(\iota_w\bar{\gamma})\wedge(\iota_w\bar{\gamma})=0$.
Obviously  $\gamma_0=(1/2)(\gamma+\bar{\gamma})$. Let $v\in V$ be such that
$(\iota_v\gamma_0)\wedge(\iota_v\gamma_0)=0$. Then
$$
0=(\iota_v\gamma_0)\wedge(\iota_v\gamma_0)=\frac{1}{4}(\iota_v\gamma+
\iota_v\bar{\gamma})\wedge(\iota_v\gamma+\iota_v\bar{\gamma})=
\frac{1}{2}(\iota_v\gamma)\wedge(\iota_v\bar{\gamma}).
$$
But $\iota_v\gamma$ is a~form of type $(2,0)$ and $\iota_v\bar{\gamma}$
a form of type $(0,2)$. Consequently the last wedge product vanishes if
and only if either $\iota_v\gamma=0$ or $\iota_v\bar{\gamma}=0$. By
virtue of the preceding lemma this implies that $v=0$.
\enddemo

Lemma 2 shows that
the both forms $\gamma_0|V$ and
$\gamma_1|V$ are of type 2.
As a~final result of
the above considerations we get the following result.
\proclaim{3. Corollary}
Let $\gamma$ be a~3-form on $V^\C$ of the type $(3,0)$. Then the
real 3-forms $(\Re\gamma)|V$ and $(\Im\gamma)|V$ on $V$ are multisymplectic
and of type 2.
\endproclaim

Let $\omega$ be a~3-form on $V$ such that $\Delta(\omega)=\{0\}$. This
means that for any $v\in V$, $v\ne0$ there is $(\iota_v\omega)\wedge
(\iota_v\omega)\ne0$. This implies that $\rank\iota_v\om\geq4$. On the other
hand obviously $\rank\iota_v\om\leq4$. Consequently, for any $v\ne0$
$\rank\iota_v\om=4$. Thus
the kernel $K(\iota_v\om)$ of the
2-form $\iota_v\omega$ has dimension 2. Moreover $v\in K(\iota_v\om)$.
Now we fix a~non-zero 6-form on $\theta$ on $V$. For any $v\in V$
there exists a~unique vector $Q(v)\in V$ such that
$$
(\iota_v\omega)\wedge\omega=\iota_{Q(v)}\theta.
$$
The mapping $Q:V\rightarrow V$ is obviously a~homomorphism. If
$v\ne0$ then $(\iota_v\omega)\wedge\omega\ne0$, and
$Q$ is an automorphism. It is also obvious that if
$v\ne0$, then the vectors $v$ and $Q(v)$ are linearly independent (apply
$\iota_v$ to the last equality). We evaluate $\iota_{Q(v)}$ on the last equality
and we get
$$
\eqalign{
(\iota_{Q(v)}\iota_v\om)\wedge\omega+(\iota_v\omega)\wedge(\iota_{Q(v)}
\omega)&=0\cr
-(\iota_v\iota_{Q(v)}\om)\wedge\omega+(\iota_v\omega)\wedge(\iota_{Q(v)}
\omega)&=0\cr
-\iota_v[(\iota_{Q(v)}\omega)\wedge\omega]+2(\iota_v\omega)\wedge
(\iota_{Q(v)}\omega)&=0}
$$
Now, apply $\iota_v$ to the last equality:
$$
(\iota_v\omega)\wedge(\iota_v\iota_{Q(v)}\omega)=0.
$$
If the 1-form $\iota_v\iota_{Q(v)}\omega$ were not the zero one
then  it would exist a~1-form $\sigma$ such that $\iota_v\omega=
\sigma\wedge\iota_v\iota_{Q(v)}\omega$, and we would get
$$
(\iota_v\omega)\wedge(\iota_v\omega)=\sigma\wedge\iota_v\iota_{Q(v)}\omega
\wedge\sigma\wedge\iota_v\iota_{Q(v)}\omega=0,
$$
which is a~contradiction. Thus we have proved the following lemma.
\proclaim{4. Lemma}
For any $v\in V$ there is $\iota_{Q(v)}\iota_v\omega=0$, i\. e\. $Q(v)\in
K(\iota_v\omega)$.
\endproclaim
This lemma shows that if $v\ne0$, then $K(\iota_v\omega)=[v,Q(v)]$. Applying
$\iota_{Q(v)}$ to the equality $(\iota_v\omega)\wedge\omega=\iota_{Q(v)}
\theta$ and using the last lemma we obtain easily the following result.
\proclaim{5. Lemma}
For any $v\in V$ there is $(\iota_v\omega)\wedge(\iota_{Q(v)}\omega)=0$.
\endproclaim

Lemma 4 shows that $v\in K(\iota_{Q(v)}\omega)$. Because $v$ and $Q(v)$ are
linearly independent, we can see that
$$
K(\iota_{Q(v)}\omega)=[v,Q(v)]=K(\iota_v\omega).
$$
If $v\ne0$, then $Q^2(v)\in K(\iota_{Q(v)}\omega)$, and consequently
there are $a(v),b(v)\in\Bbb R$ such that
$$
Q^2(v)=a(v)v+b(v)Q(v).
$$
For any $v\in V$
$$
(\iota_{Q(v)}\omega)\wedge\omega=\iota_{Q^2(v)}\theta.
$$
Let us assume that $v\ne0$. Then
$$
(\iota_{Q(v)}\omega)\wedge\omega=a(v)\iota_v\theta+b(v)\iota_{Q(v)}\theta,
$$
and applying $\iota_v$ we obtain $b(v)\iota_v\iota_{Q(v)}\theta=0$, which shows
that $b(v)=0$ for any $v\ne0$. Consequently,  $Q^2(v)=a(v)v$ for any
$v\ne0$.
\proclaim{6. Lemma}
Let $A:V\rightarrow V$ be an automorphism, and $a:V\backslash\{0\}\rightarrow
\Bbb R$ a~function such that
$$
A(v)=a(v)v\quad\text{for any }v\ne0.
$$
Then the function $a$ is constant.
\endproclaim
\demo{Proof}
The condition on $A$ means that every vector $v$ of $V$ is an eigenvector
of $A$ with the eigenvalue $a(v)$. But the eigenvalues of two different
vectors have to be the same otherwise their sum would not be an
eigenvector.
\enddemo

Applying Lemma 6 on $Q^2$ we get $Q^2=aI$. If $a>0$, then $V=V^+\oplus V^-$, and
$$
Qv=\sqrt{a}v\text{ for }v\in V^+,\quad Qv=-\sqrt{a}v\text{ for }v\in
V^-.
$$
At least one of the subspaces $V^+$ and $V^-$ is non-trivial. Let us
assume for example that $V^+\ne\{0\}$. Then there is $v\in V^+$,
$v\ne0$, and $Qv=\sqrt{a}v$, which is a~contradiction because
the vectors $v$ and $Qv$ are linearly independent. This proves
that $a<0$. We can now see that the automorphisms
$$
J_+=\frac{1}{\sqrt{-a}}Q\text{ and }J_-=-\frac{1}{\sqrt{-a}}Q
\text{ satisfy }J_+^2=-I\text{ and }J_-^2=-I,
$$
i\. e\. they define complex structures on $V$ and $J_-=-J_+$.
Setting
$$
\theta_+=\sqrt{-a}\theta,\quad\theta_-=-\sqrt{-a}\theta
$$
we get
$$
(\iota_v\omega)\wedge\omega=\iota_{J_+v}\theta_+,\quad
(\iota_v\omega)\wedge\omega=\iota_{J_-v}\theta_-.
$$

In the sequel we shall denote $J=J_+$. The same results which are valid
for $J_+$ hold also for $J_-$.
\proclaim{7. Lemma}
There exists  a~unique (up to the the sign) complex structure $J$ on $V$
such that the form $\omega$ satisfies the relation
$$
\omega(Jv_1,v_2,v_3)=\omega(v_1,Jv_2,v_3)=\omega(v_1,v_2,Jv_3)\quad
\text{for any }v_1,v_2,v_3\in V.
$$
We recall that such a form $\omega$ is usually called pure with respect to
$J$.

\endproclaim
\demo{Proof}
We shall prove first that the complex structure $J$ defined above
satisfies the relation. By virtue of Lemma 4 for any $v,v'\in V$
$\omega(v,Jv,v')=0$. Therefore we get
$$
\split
0=\omega(v_1+v_2,J(v_1+v_2),v_3)&=\omega(v_1,Jv_2,v_3)+\omega(v_2,Jv_1,v_3)\\
&=-\omega(Jv_1,v_2,v_3)+\omega(v_1,Jv_2,v_3),\\
\endsplit
$$
which gives
$$
\omega(Jv_1,v_2,v_3)=\omega(v_1,Jv_2,v_3).
$$
Obviously, the opposite complex structure $-J$ satisfies the same
relation. We prove that there is no other complex
structure with the same property. Let $\tilde{J}$ be a~complex structure
on $V$ satisfying the above relation. We set $A=\tilde{J}J^{-1}$. Then
we get
$$
\split
\omega(v_1,Av_2,Av_3)=\omega(v_1,\tilde{J}Jv_2,\tilde{J}Jv_3)=&
\omega(v_1,Jv_2,\tilde{J}^2Jv_3)=-\omega(v_1,Jv_2,Jv_3)\\
=&-\omega(v_1,v_2,J^2v_3)=\omega(v_1,v_2,v_3).
\endsplit
$$
Any automorphism $A$ satisfying this
identity is $\pm I$. Really, the identity means that that $A$ is an
automorphism of the 2-form $\iota_v\omega$. Consequently, $A$ preserves
the kernel $K(\iota_v\omega)=[v,Jv]$. On the other hand it is obvious
that any subspace of the form $[v,Jv]$ is the kernel of $\iota_v\omega$.
Considering $V$ as a~complex vector space with the complex structure
$J$, we can say that every 1-dimensional complex subspace is the kernel
of the 2-form $\iota_v\omega$ for some $v\in V$, $v\ne0$, and
consequently is invariant under the automorphism $A$.
Similarly as in Lemma 6 we conclude, that $A=\lambda I$, $\lambda\in\C$.
If we write $\lambda=\lambda_0
+i\lambda_1$, then $A=\lambda_0I+\lambda_1J$ and
$$
\split
&\omega(v_1,v_2,v_3)=\omega(v_1,Av_2,Av_3)\\
=&\omega(v_1,\lambda_0v_2+\lambda_1Jv_2,\lambda_0v_3+\lambda_1Jv_3)\\
=&\lambda_0^2\omega(v_1,v_2,v_3)+\lambda_0\lambda_1\omega(v_1,v_2,Jv_3)+
\lambda_0\lambda_1\omega(v_1,Jv_2,v_3)+\lambda_1^2\omega(v_1,Jv_2,Jv_3)\\
\endsplit
$$
$$
(\lambda_0^2-\lambda_1^2-1)\omega(v_1,v_2,v_3)+
2\lambda_0\lambda_1\omega(v_1,v_2,Jv_3)=0
$$
We shall use this last equation together with another one obtained by
writing $Jv_3$ instead of $v_3$. In this way we get the
system
$$
\split
(\lambda_0^2-\lambda_1^2-1)\omega(v_1,v_2,v_3)+
2\lambda_0\lambda_1\omega(v_1,v_2,Jv_3)&=0\\
-2\lambda_0\lambda_1\omega(v_1,v_2,v_3)+
(\lambda_0^2-\lambda_1^2-1)\omega(v_1,v_2,Jv_3)&=0.\\
\endsplit
$$
Because it has a~non-trivial solution there must be
$$
\vmatrix
\lambda_0^2-\lambda_1^2-1 & 2\lambda_0\lambda_1\\
-2\lambda_0\lambda_1       & \lambda_0^2-\lambda_1^2-1
\endvmatrix
=0.
$$
It is easy to verify that the solution of the last equation is
$\lambda_0=\pm1$ and $\lambda_1=0$. This finishes the proof.
\enddemo

We shall now consider the vector space $V$ together with a~complex
structure $J$, and a~3-form $\omega$ on $V$ which is pure with respect to
this complex structure. Firstly we define a~real 3-form
$\gamma_0$ on $V^\C$. We set
\roster
\item"{}"$\gamma_0(v_1,v_2,v_3)=\omega(v_1,v_2,v_3)$,
\item"{}"$\gamma_0(iv_1,v_2,v_3)=\omega(Jv_1,v_2,v_3)$,
\item"{}"$\gamma_0(iv_1,iv_2,v_3)=\omega(Jv_1,Jv_2,v_3)$,
\item"{}"$\gamma_0(iv_1,iv_2,iv_3)=\omega(Jv_1,Jv_2,Jv_3)$,
\endroster
for $v_1,v_2,v_3\in V$. Then $\gamma_0$ extends uniquely to
a real 3-form on $V^\C$.
We can find easily that
$$
\gamma_0(iw_1,w_2,w_3)=\gamma_0(w_1,iw_2,w_3)=\gamma_0(w_1,w_2,iw_3)
$$
for any $w_1,w_2,w_3\in V^\C$. Further, we set
$$
\gamma_1(w_1,w_2,w_3)=-\gamma_0(iw_1,w_2,w_3)\quad\text{for }w_1,w_2,w_3
\in V^\C.
$$
It is obvious that $\gamma_1$ is a~real 3-form satisfying
$$
\gamma_1(iw_1,w_2,w_3)=\gamma_1(w_1,iw_2,w_3)=\gamma_1(w_1,w_2,iw_3)
$$
for any $w_1,w_2,w_3\in V^\C$. Now we define
$$
\gamma(w_1,w_2,w_3)=\gamma_0(w_1,w_2,w_3)+i\gamma_1(w_1,w_2,w_3)
\quad\text{for }w_1,w_2,w_3\in V^\C.
$$
It is obvious that $\gamma$ is skew symmetric and 3-linear over $\Bbb R$
and has complex values. Moreover
$$
\split
\gamma(iw_1,w_2,w_3)&=\gamma_0(iw_1,w_2,w_3)+i\gamma_1(iw_1,w_2,w_3)\\
=-\gamma_1(w_1,w_2,w_3)&-i\gamma_0(i^2w_1,w_2,w_3)=
-\gamma_1(w_1,w_2,w_3)+i\gamma_0(w_1,w_2,w_3)=\\
&=i[\gamma_0(w_1,w_2,w_3)+i\gamma_1(w_1,w_2,w_3)]=i\gamma(w_1,w_2,w_3),\\
\endsplit
$$
which proves that $\gamma$ is a~complex 3-form on $V^\C$. Now we prove
that $\gamma$ is a~form of type $(3,0)$. Obviously, it suffices to
prove that for $v_1+iJv_1\in V^{0,1}$ and $v_2,v_3\in V$ there is
$\gamma(v_1+iJv_1,v_2,v_3)=0$. Really,
$$
\split
&\gamma(v_1+iJv_1,v_2,v_3)=\gamma(v_1,v_2,v_3)+i\gamma(Jv_1,v_2,v_3)\\
=&\gamma_0(v_1,v_2,v_3)+i\gamma_1(v_1,v_2,v_3)+i\gamma_0(Jv_1,v_2,v_3)-
\gamma_1(Jv_1,v_2,v_3)\\
=&\gamma_0(v_1,v_2,v_3)-i\gamma_0(iv_1,v_2,v_3)+i\gamma_0(Jv_1,v_2,v_3)
+\gamma_0(iJv_1,v_2,v_3)].\\
\endsplit
$$
Now $\gamma_0(iJv_1,v_2,v_3)]=\om(J^2v_1,v_2,v_3)=-\om(v_1,v_2,v_3)=
-\gamma_0(v_1,v_2,v_3)$ and the real part of the last expression is zero,
further $\gamma_0(Jv_1,v_2,v_3)=\om(Jv_1,v_2,v_3)=\gamma_0(iv_1,v_2,v_3)$
and the complex part of the expression is zero as well.
Now we get easily the following proposition.

\proclaim{8. Proposition}
Let $\omega$ be a~real 3-form on $V$ satisfying $\Delta(\omega)=\{0\}$,
and let $J$ be a~complex structure on $V$ (one of the two) such that
$$
\omega(Jv_1,v_2,v_3)=\omega(v_1,Jv_2,v_3)=\omega(v_1,v_2,Jv_3).
$$
Then there exists on $V^\C$ a~unique complex 3-form $\gamma$ of type $(3,0)$
such that
$$
\omega=(\Re\gamma)|V.
$$
\endproclaim

\remark{Remark} The complex structure $J$ on $V$ can be introduced also by means of
the Hitchin's invariant $\lambda$, as in [H]. Forms of type 2
form an open subset $U$ in $\Lambda^3V^*$.
Hitchin
has shown that this manifold also carries an almost complex structure, which is
integrable.
Hitchin uses the following way to introduce an almost complex structure on
$U$. 
$U\subset\Lambda^3V^*$ can be seen as a symplectic manifold
(let $\theta$ be a fixed element in $\Lambda^6V^*$;
one defines the symplectic form $\Theta$ on $\Lambda^3V^*$
by the equation $\w_1\wedge\w_2=\Theta(\w_1,\w_2)\theta$).
Then the derivative of the Hamiltonian vector
field corresponding to the function $\sqrt{-\lambda(\w)}$ on $U$ gives
an integrable almost complex structure on $U$.
That was for the Hitchin's construction.

There is another way of introducing the (Hitchin's) almost complex structure
on $U$. Given a $3$-form $\w\in U$ we choose the complex structure $J_\w$ on
$V$ (one of the two), whose existence is guaranteed by the lemma 7.
Then we define
 endomorphisms $A_{J_\w}$ and $D_{J_\w}$ of $\Lambda^kV^*$
by
$$
(A_{J_\w}{\Omega})(v_1,\dots,v_k)=\Omega(J_\w v_1,\dots,J_\w v_k),
$$
$$
(D_{J_\w}{\Omega})(v_1,\dots,v_k)=\sum_{i=1}^k\Omega(v_1,\dots,v_{i-1},J_\w v_i,v_{i+1},
\dots,v_k).
$$
Then $A_{J_\w}$ is an automorphism of $\Lambda V^*$  and $D_{J_\w}$ is a derivation
of $\Lambda V^*$. If $k=3$  
then the automorphism $-\frac12(A_{J_\w}+D_{J_\w})$ of $\Lambda^3V^*$
($=T_\w U$)
gives a complex structure on $U$
and coincides with the Hitchin's one.

\head
\centerline{\bf $3$-forms on manifolds}
\endhead
We use facts from the previous section to obtain some global results
on $3$-forms on $6$-dimensional manifolds.
We shall denote by $X$, $Y$, $Z$ the real vector fields on a (real) manifold $M$ and
by $V$, $W$ the complex vector fields on $M$. $\Cal X(M)$ stands for the set
of all (real) vector fields on $M$, $\Cal X^\C(M)$ means all the complex
vector fields on $M$.

A 3-form
$\omega$ on $M$ is called the form of type 2 if for every $x\in M$ there is
$\Delta(\omega_x)=\{0\}$.
Let $\w$ be a form of type 2 on $M$ and let $U\subset M$ be an open orientable
submanifold. Then there exists an everywhere nonzero
differentiable 6-form on $U$.  In each $T_xM$, $x\in U$
construct $J_-$
and $J_+$ as in Lemma 7. The construction is evidently smooth on $U$. Thus

\proclaim{9. Lemma}
Let $\om$ be a~form of type 2 on $M$ and
 let $U\subset M$ be an orientable open submanifold.
Then there exist two differentiable almost complex structures $J_+$ and $J_-$ on $U$ such
that
\roster
\item"{(i)}"  $J_++J_-=0$,
\item"{(ii)}" $\omega(J_+X_1,X_2,X_3)=\omega(X_1,J_+X_2,X_3)=
               \omega(X_1,X_2,J_+X_3)$,
\item"{(iii)}" $\omega(J_-X_1,X_2,X_3)=\omega(X_1,J_-X_2,X_3)=
               \omega(X_1,X_2,J_-X_3)$,
\endroster
for any vector fields $X_1$, $X_2$, $X_3$.
\endproclaim

At each point $x\in M$ consider a~1-dimensional subspace of the
space $\Cal T_{1x}^1(M)$ of tensors of type $(1,1)$ at $x$ generated by the
tensors $J_{+x}$ and $J_{-x}$. The above considerations show that
it is  a~1-dimensional subbundle $\Cal J\subset\Cal
T_1^1(M)$.
\proclaim{10. Lemma}
The 1-dimensional vector bundles $\Cal J$ and $\Lambda^6T^*(M)$ are
isomorphic.
\endproclaim
\demo{Proof}
Let us choose a~riemannian metric $g_0$ on $TM$. If $x\in M$ and
$v,v'\in T_xM$ we define a~riemannian metric $g$ by the formula
$$
g(v,v')=g_0(v,v')+g_0(J_+v,J_+v')=g_0(v,v')+g_0(J_-v,J_-v').
$$
It is obvious that for any $v,v'\in T_xM$ we have
$$
g(J_+v,J_+v')=g(v,v'),\quad g(J_-v,J_-v')=g(v,v').
$$
We now define
$$
\sigma_{+}(v,v')=g(J_+v,v'),\quad \sigma_{-}(v,v')=g(J_-v,v').
$$
It is easy to verify that $\sigma_{+}$ and $\sigma_{-}$ are nonzero
2-forms on $T_xM$ satisfying $\sigma_++\sigma_-=0$.

We define an isomorphism $h:\Cal J\rightarrow\Lambda^6
T^*M$. Let $x\in M$ and let $A\in\Cal J_x$. We can write
$$
A=aJ_+,\quad A=-aJ_-.
$$
We set
$$
hA=a\sigma_+\wedge\sigma_+\wedge\sigma_+=
-a\sigma_-\wedge\sigma_-\wedge\sigma_-.
$$

\proclaim{11. Corollary}
There exist two almost complex structures $J_+$ and $J_-$ on $M$ such
that
\roster
\item"{(i)}"  $J_++J_-=0$,
\item"{(ii)}" $\omega(J_+X_1,X_2,X_3)=\omega(X_1,J_+X_2,X_3)=
               \omega(X_1,X_2,J_+X_3)$,
\item"{(iii)}" $\omega(J_-X_1,X_2,X_3)=\omega(X_1,J_-X_2,X_3)=
               \omega(X_1,X_2,J_-X_3)$,
\endroster
for any vector fields $X_1$, $X_2$, $X_3$ if and only if the manifold $M$ is orientable.
\endproclaim

Hence the assertions in the rest of the article can be simplified correspondingly
if $M$ is an orientable manifold.

\proclaim{12. Lemma}
Let $J$ be an almost complex structure on $M$ such that for any vector
fields $X_1,X_2,X_3\in\Cal X(M)$ there is
$$
\omega(JX_1,X_2,X_3)=\omega(X_1,JX_2,X_3)=\omega(X_1,X_2,JX_3).
$$
If $\nabla$ is a~linear connection on $M$ such that $\nabla\omega=0$,
then also $\nabla J=0$.
\endproclaim
\demo{Proof}
Let $Y\in\Cal X(M)$, and let us consider the covariant derivative
$\nabla_Y$. We get
$$
\split
0=(\nabla_Y&\omega)(JX_1,X_2,X_3)=Y(\omega(JX_1,X_2,X_3)
-\omega((\nabla_YJ)X_1,X_2,X_3)\\
-&\omega(J\nabla_YX_1,X_2,X_2)
-\omega(JX_1,\nabla_YX_2,X_3)-\omega(JX_1,X_2,\nabla_YX_3),\\
0=(\nabla_Y&\omega)(X_1,JX_2,X_3)=Y(\omega(JX_1,X_2,X_3)
-\omega(\nabla_YX_1,JX_2,X_3)\\
-&\omega(X_1,(\nabla_YJ)X_2,X_3)
-\omega(X_1,J\nabla_YX_2,X_3)-\omega(X_1,JX_2,\nabla_YX_3).\\
\endsplit
$$
Because the above expressions are equal we find easily that
$$
\omega((\nabla_YJ)X_1,X_2,X_3)=\omega(X_1,(\nabla_YJ)X_2,X_3).
$$
We denote $A=\nabla_YJ$. Extending in the obvious way the above equality,
we get
$$
\omega(AX_1,X_2,X_3)=\omega(X_1,AX_2,X_3)=\omega(X_1,X_2,AX_3).
$$
Moreover $J^2=-I$, and applying $\nabla_Y$ to this equality, we get
$$
AJ+JA=0.
$$
We know that $K(\iota_X\omega)=[X,JX]$. Furthermore
$$
\omega(X,AX,X')=\omega(X,X,AX')=0,\quad
\omega(X,AJX,X')=\omega(X,JX,AX')=0,
$$
which shows that $A$ preserves the distribution $[X,JX]$. By the very same
arguments as in
Lemma 7 we can see that $A=\lambda_0I+\lambda_1J$. Consequently
$$
\gather
(\lambda_0I+\lambda_1J)J+J(\lambda_0I+\lambda_1J)=0\\
-2\lambda_1I+2\lambda_0J=0,
\endgather
$$
which implies $\lambda_0=\lambda_1=0$. Thus
$\nabla_YJ=A=0$.
\enddemo

The statement of the previous lemma can be in a way reversed,
and we get

\proclaim{13. Proposition}
Let $\omega$ be a~real 3-form on a~6-dimensional differentiable manifold
$M$ satisfying $\Delta(\omega_x)=\{0\}$ for any $x\in M$.
Let $J$
be an almost complex structure on $M$ such that for any vector fields
$X_1,X_2,X_3\in\Cal X(M)$ there is
$$
\omega(JX_1,X_2,X_3)=\omega(X_1,JX_2,X_3)=\omega(X_1,X_2,JX_3).
$$
Then there exists a~symmetric connection $\tilde{\nabla}$ on $M$ such that
$\tilde{\nabla}\omega=0$ if and only if the following conditions are
satisfied
\roster
\item"{(i)}" $d\omega=0$,
\item"{(ii)}" the almost complex structure $J$ is integrable.
\endroster
\endproclaim
\demo{Proof}
First, we prove that the integrability of the structure $J$ and the fact that
$\omega$ is closed implies the existence of a symmetric connection with
respect to which $\omega$ is parallel.

For any connection $\nabla$ on $M$ we shall denote by the same symbol
its complexification. Namely, we set
$$
\nabla_{X_0+iX_1}(Y_0+iY_1)=(\nabla_{X_0}Y_0-\nabla_{X_1}Y_1)+
i(\nabla_{X_0}Y_1+\nabla_{X_1}Y_0).
$$
Let us assume that there exists a~symmetric connection $\overset{\circ}
\to\nabla$ such that $\overset{\circ}\to\nabla J=0$. We shall consider
a 3-form $\gamma$ of type $(3,0)$ such that $(\Re\gamma)|TM=\omega$. Our
next aim is to try to find a~symmetric connection
$$
\nabla_VW=\overset{\circ}\to\nabla_VW+Q(V,W)
$$
satisfying $\nabla_V\gamma=0$. Obviously, the connection $\nabla$ is
symmetric if and only if
$$
Q(V,W)=Q(W,V).
$$
Moreover, $\nabla_V\gamma=0$ hints that
$\nabla J=0$.
$$
0=(\nabla_VJ)W=\nabla_V(JW)-J\nabla_VW=\overset{\circ}\to\nabla_V(JW)+
Q(V,JW)-J\overset{\circ}\to\nabla_VW-JQ(V,W),
$$
which shows that we should require
$$
Q(JV,W)=Q(V,JW)=JQ(V,W).
$$
Because $\overset{\circ}\to\nabla J=0$, we can immediately see that
for any $V\in\Cal X^\C(M)$ the covariant derivative $\overset{\circ}\to
\nabla_V\gamma$ is again a~form of type $(3,0)$. Consequently there
exists a~uniquely determined complex 1-form $\rho$ such that
$$
\overset{\circ}\to\nabla_V\gamma=\rho(V)\gamma.
$$
Then
$$
\split
&(\nabla_V\gamma)(W_1,W_2,W_3)\\
=&V(\gamma(W_1,W_2,W_3))
-\gamma(\nabla_VW_1,W_2,W_3)-\gamma(W_1,\nabla_VW_2,W_3)
-\gamma(W_1,W_2,\nabla_VW_3)\\
=&V(\gamma(W_1,W_2,W_3))-\gamma(\overset{\circ}\to\nabla_VW_1,W_2,W_3)
-\gamma(W_1,\overset{\circ}\to\nabla_VW_2,W_3)
-\gamma(W_1,W_2,\overset{\circ}\to\nabla_VW_3)\\
&-\gamma(Q(V,W_1),W_2,W_3)-\gamma(W_1,Q(V,W_2),W_3)
-\gamma(W_1,W_2,Q(V,W_3))\\
=&\rho(V)\gamma(W_1,W_2,W_3)\\
&-\gamma(Q(V,W_1),W_2,W_3)-\gamma(W_1,Q(V,W_2),W_3)
-\gamma(W_1,W_2,Q(V,W_3)).
\endsplit
$$
In other words $\nabla_V\gamma=0$ if and only if
$$
\split
&\rho(V)\gamma(W_1,W_2,W_3)\\
=\gamma(Q(V,W_1),W_2,W_3)+&\gamma(W_1,Q(V,W_2),W_3)
+\gamma(W_1,W_2,Q(V,W_3)).
\endsplit
$$

\proclaim{Sublemma}
If $d\gamma=0$, then $\rho$ is a~form of type $(1,0)$.
\endproclaim
\demo{Proof}
Let $V_1\in T^{0,1}(M)$. Because $\overset{\circ}\to\nabla$ is symmetric
$d\gamma=
-\Cal A(\overset{\circ}\to\nabla\gamma)$, where $\Cal A$ denotes the
alternation. We obtain
$$
\split
0=-4!(d\gamma)(V_1,V_2,V_3,V_4)&=\sum_{\pi}\sign(\pi)(\overset{\circ}\to
\nabla_{V_{\pi1}}\gamma)(V_{\pi2},V_{\pi3},V_{\pi4})\\
+\sum_\tau\sign(\tau)(\overset{\circ}\to\nabla_{V_1}\gamma)(V_{\tau2},
V_{\tau3},V_{\tau4})&=3!(\overset{\circ}\to\nabla_{V_1}\gamma)(V_2,V_3,V_4)=
3!\rho(V_1)\gamma(V_2,V_3,V_4).
\endsplit
$$
The first sum is taken over all permutations $\pi$ satisfying
$\pi1>1$, and the second one is taken over all permutations of the
set $\{2,3,4\}$. The first sum obviously vanishes, and
$\rho(V_1)=0$. This finishes the proof.
\enddemo

We set now
$$
Q(V,W)=\frac{1}{8}[\rho(V)W-\rho(JV)JW+\rho(W)V-\rho(JW)JV].
$$
It is easy to see that $Q(JV,W)=Q(V,JW)=JQ(V,W)$. For $V,W_1,
W_2,W_3\in T^{1,0}(M)$ we can compute
$$
\split
&8\gamma(Q(V,W_1),W_2,W_3)\\
=&\gamma(\rho(V)W_1-\rho(JV)JW_1+\rho(W_1)V-\rho(JW_1)JV,W_2,W_3)\\
=&\gamma(2\rho(V)W_1+2\rho(W_1)V,W_2,W_3)=2\rho(V)\gamma(W_1,W_2,W_3)+
2\rho(W_1)\gamma(V,W_2,W_3),\\
\endsplit
$$
where we used for $V\in T^{(1,0)}(M)$ that $\rho(JV)=i\rho(V)$ and $\gamma(JV,V',V'')=i\gamma(V,V',V'')$,
since $\gamma$ is of type $(3,0)$ and $\rho$
of type $(1,0)$.

Similarly we can compute $\gamma(W_1,Q(V,W_2),W_3)$ and $\gamma(W_1,
W_2,Q(V,W_3))$. Without a~loss of generality we can assume that the
vector fields $W_1,W_2,W_3$ are linearly independent (over $\Bbb C$).
Then we can find uniquely determined complex functions $f_1,f_2,f_3$
such that
$$
V=f_1W_1+f_2W_2+f_3W_3.
$$
Then we get
$$
\split
&\rho(W_1)\gamma(V,W_2,W_3)+\rho(W_2)\gamma(W_1,V,W_3)+
\rho(W_3)\gamma(W_1,W_2,V)\\
=&f_1\rho(W_1)\gamma(W_1,W_2,W_3)+f_2\rho(W_2)\gamma(W_1,W_2,W_3)+
f_3\rho(W_3)\gamma(W_1,W_2,W_3)\\
=&\rho(f_1W_1+f_2W_2+f_3W_3)\gamma(W_1,W_2,W_3)=
\rho(V)\gamma(W_1,W_2,W_3).
\endsplit
$$
Finally we obtain
$$
\split
\gamma(Q(V,W_1),W_2,W_3)+&\gamma(W_1,Q(V,W_2),W_3)+
\gamma(W_1,W_2,Q(V,W_3))\\
=&\rho(V)\gamma(W_1,W_2,W_3).\\
\endsplit
$$
which proves $\nabla_V\gamma=0$.

Let us continue in the main stream of the proof. We shall now use the  complex connection
$\nabla$.
For $X,Y\in TM$ we
shall denote $\nabla^0_XY=\Re\nabla_XY$ and $\nabla^1_XY=\Im\nabla_XY$.
This means that we have $\nabla_XY=\nabla^0_XY+i\nabla^1_XY$. For a real
function $f$ on $M$ we have
$$
\nabla_X(fY)=\nabla^0_X(fY)+i\nabla^1_X(fX),
$$
$$
\nabla_X(fY)=(Xf)Y+f\nabla_XY=[(Xf)Y+f\nabla^0_XY]+if\nabla^1_XY,
$$
which implies
$$
\nabla^0_X(fY)=(Xf)Y+f\nabla^0_XY,\quad
\nabla^1_X(fY)=f\nabla^1_XY.
$$
This shows that $\nabla^0$ is a real connection while $\nabla^1$ is a
real tensor field of type $(1,2)$. We have also
$$
0=\nabla_XY-\nabla_YX-[X,Y]=\nabla^0_XY+i\nabla^1_XY-\nabla^0_YX
-i\nabla^1_YX-[X,Y]=
$$
$$
=[\nabla^0_XY-\nabla^0_YX-[X,Y]]+i[\nabla^1_XY-\nabla^1_YX],
$$
which shows that
$$
\nabla^0_XY-\nabla^0_YX-[X,Y]=0,\quad \nabla^1_XY-\nabla^1_YX=0.
$$
These equations show that the connection $\nabla^0$ is symmetric, and
that the tensor $\nabla^1$ is also symmetric. Moreover, we have
$$
\nabla_X(JY)=\nabla^0_X(JY)+i\nabla^1_X(JY),
$$
$$
\nabla_X(JY)=J\nabla_XY=J\nabla^0_XY+iJ\nabla^1_XY,
$$
which gives
$$
\nabla^0_XJ=0,\quad \nabla^1_X(JY)=J\nabla^1_XY.
$$

For the real vectors $X,Y_1,Y_2,Y_3\in TM$ we can compute
$$
0=(\nabla_X\gamma)(Y_1,Y_2,Y_3)=X(\gamma(Y_1,Y_2,Y_3))-
$$
$$
-\gamma(\nabla_XY_1,Y_2,Y_3)-\gamma(Y_1,\nabla_XY_2,Y_3)
-\gamma(Y_1,Y_2,\nabla_XY_3)=
$$
$$
=X(\gamma(Y_1,Y_2,Y_3))-
$$
$$
-\gamma(\nabla_X^0Y_1+i\nabla_X^1Y_1,Y_2,Y_3)
-\gamma(Y_1,\nabla^0_XY_2+i\nabla^1_XY_2,Y_3)
-\gamma(Y_1,Y_2,\nabla^0_XY_3+i\nabla^1_XY_3)=
$$
$$
=X(\gamma(Y_1,Y_2,Y_3))-\gamma(\nabla_X^0Y_1,Y_2,Y_3)
-\gamma(Y_1,\nabla^0_XY_2,Y_3)-\gamma(Y_1,Y_2,\nabla^0_XY_3)-
$$
$$
-i[\gamma(\nabla_X^1Y_1,Y_2,Y_3)
+\gamma(Y_1,\nabla^1_XY_2,Y_3)+\gamma(Y_1,Y_2,\nabla^1_XY_3)]=
$$
$$
=[X(\gamma_0(Y_1,Y_2,Y_3))-\gamma_0(\nabla_X^0Y_1,Y_2,Y_3)
-\gamma_0(Y_1,\nabla^0_XY_2,Y_3)-\gamma_0(Y_1,Y_2,\nabla^0_XY_3)+
$$
$$
+\gamma_1(\nabla_X^1Y_1,Y_2,Y_3)
+\gamma_1(Y_1,\nabla^1_XY_2,Y_3)+\gamma_1(Y_1,Y_2,\nabla^1_XY_3)]+
$$
$$
+i[X(\gamma_1(Y_1,Y_2,Y_3))-\gamma_1(\nabla_X^0Y_1,Y_2,Y_3)
-\gamma_1(Y_1,\nabla^0_XY_2,Y_3)-\gamma_1(Y_1,Y_2,\nabla^0_XY_3)-
$$
$$
-\gamma_0(\nabla_X^1Y_1,Y_2,Y_3)-\gamma_0(Y_1,\nabla^1_XY_2,Y_3)
-\gamma_0(Y_1,Y_2,\nabla^1_XY_3)].
$$
This shows that the real part (as well as the complex one, which gives
in fact the same identity) is zero. Using the relations between
$\gamma_0$ and $\gamma_1$ we get
$$
0=X(\gamma_0(Y_1,Y_2,Y_3))-\gamma_0(\nabla_X^0Y_1,Y_2,Y_3)
-\gamma_0(Y_1,\nabla^0_XY_2,Y_3)-\gamma_0(Y_1,Y_2,\nabla^0_XY_3)-
$$
$$
-\gamma_0(J\nabla_X^1Y_1,Y_2,Y_3)-\gamma_0(Y_1,J\nabla^1_XY_2,Y_3)
-\gamma_0(Y_1,Y_2,J\nabla^1_XY_3)=
$$
$$
=X(\gamma_0(Y_1,Y_2,Y_3))-\gamma_0(\nabla_X^0Y_1+J\nabla_X^1Y_1,Y_2,Y_3)
$$
$$
-\gamma_0(Y_1,\nabla^0_XY_2+J\nabla^1_XY_2,Y_3)
-\gamma_0(Y_1,Y_2,\nabla^0_XY_3+J\nabla^1_XY_3).
$$

We define now
$$
\tilde{\nabla}_XY=\nabla^0_XY+J\nabla^1_XY.
$$
It is easy to verify that $\tilde{\nabla}$ is a real connection. Moreover,
the previous equation shows that
$$
\tilde{\nabla}\gamma_0=0.
$$
Furthermore, it is very easy to see that the connection
$\tilde{\nabla}$ is symmetric.

The inverse implication can be proved easily.
\enddemo

Let us use the standard definition of integrability of a~$k$-form $\omega$ on $M$,
that is every $x\in M$ has a~neighbourhood $N$ such that $\omega$ has the
constant expresion in $dx^i$, $x^i$ being suitable coordinate functions on
$N$.

\proclaim{14. Corollary}
Let $\omega$ be a~real 3-form on a~6-dimensional differentiable manifold
$M$ satisfying $\Delta(\omega_x)=\{0\}$ for any $x\in M$.
Let $J$
be an almost complex structure on $M$ such that for any vector fields
$X_1,X_2,X_3\in\Cal X(M)$ there is
$$
\omega(JX_1,X_2,X_3)=\omega(X_1,JX_2,X_3)=\omega(X_1,X_2,JX_3).
$$
Then $\om$
is integrable if and only if there exists a~symmetric connection $\nabla$
preserving $\om$, that is $\nabla\om=0$.
\endproclaim
\demo{Proof}
Let $\nabla$ be a~symmetric connection such that  $\nabla\om=0$. Then according
to the previous proposition $d\om=0$ and $J$ is integrable. Then we construct
the complex form $\gamma$ on $T^\C M$
of type $(3,0)$ such  and $\Re\gamma|T_xM=\om$, for any $x\in M$
(point by point, according to Proposition 8). Moreover if $\om$ is closed
than so is $\gamma$.
That is $\gamma=f\cdot dz^1\wedge dz^2\wedge dz^3$, where $z^1$,
$z^2$, and $z^3$ are (complex) coordinate functions on $M$,
$dz^1$, $dz^2$, $dz^3$ are a basis of $\Lambda^{1,0}M$ and $f$ a~function on $M$.
Further
$$
0=d\gamma=\del\gamma+\overline{\del}\gamma
=\del f\cdot dz^1\wedge dz^2\wedge dz^3
+\overline{\del} f\cdot dz^1\wedge dz^2\wedge dz^3.
$$
Evidently $\del\gamma=0$, which means $\overline{\del} f=0$ and $f$ is
holomorphic. Now we exploit a~standard trick. There exists a~holomorphic
function $F(z^1,z^2,z^3)$ such that $\frac{\del F}{\del z^1}=f$. We introduce
new complex coordinates $\tilde{z}^1=F(z^1,z^2,z^3)$,
$\tilde{z}^2=z^2$, and $\tilde{z}^3=z^3$. Then $\gamma=
fdz^1\wedge dz^2\wedge dz^3=
d\tilde{z}^1\wedge d\tilde{z}^2\wedge d\tilde{z}^3$.
Now write $\tilde{z}^1=x^1+ix^4$, $\tilde{z}^2=x^2+ix^5$, and
$\tilde{z}^3=x^3+ix^6$ for real coordinate
functions $x^1$, $x^2$, $x^3$, $x^4$, $x^5$, and $x^6$ on M. There is
$$
\split
\Re\gamma=&\Re(d(x^1+ix^4)\wedge d(x^2+ix^5)\wedge d(x^3+ix^6))\\
=&dx^1\wedge dx^2\wedge dx^3-dx^1\wedge dx^5\wedge dx^6+
dx^2\wedge dx^4\wedge dx^6-dx^3\wedge dx^4\wedge dx^5.
\endsplit
$$
And $\om=(\Re\gamma)|TM$ is an integrable  on $M$. 

Conversely, if $\om$ is integrable, then for any $x\in M$ there is a~basis
$dx_1,\dots,dx_6$ of $T^*N$ in some neighbourhoud $N\subset M$ of $x$ such
that $\om$ has constant expression in all $T_xM$, $x\in N$. Then
the flat connection $\nabla$ given by the coordinate system $x_1,\dots,x_6$
is symmetric and $\nabla\om=0$ on $N$. We use the partition of the unity
and extend $\nabla$ over the whole $M$.
\enddemo

We can reformulate the Proposition 13 as  " The Darboux theorem for
type 2 forms":

\proclaim{15. Corollary}
Let $\omega$ be a~real 3-form on a~6-dimensional differentiable manifold
$M$ satisfying $\Delta(\omega_x)=\{0\}$ for any $x\in M$.
Let $J$
be an almost complex structure on $M$ such that for any vector fields
$X_1,X_2,X_3\in\Cal X(M)$ there is
$$
\omega(JX_1,X_2,X_3)=\omega(X_1,JX_2,X_3)=\omega(X_1,X_2,JX_3).
$$
Then $\om$ is integrable
if and only if the following conditions are
satisfied
\roster
\item"{(i)}" $d\omega=0$,
\item"{(ii)}" the almost complex structure $J$ is integrable.
\endroster
\endproclaim

\remark{\bf 16. Observation}
There is an interesting relation between structures given
by a~form of type 2 on $6$-dimensional vector spaces and $G_2$-structures on
$7$-dimensional ones ($G_2$ being the exeptional Lie group, the group
of automorphisms of the algebra of Caley numbers and also the group
of automorphism of the $3$-form given below),
i.e. structures given by a~form of the type
$$
\split
&\alpha_1\wedge\alpha_2\wedge\alpha_3+
       \alpha_1\wedge\alpha_4\wedge\alpha_5-
       \alpha_1\wedge\alpha_6\wedge\alpha_7+
       \alpha_2\wedge\alpha_4\wedge\alpha_6+
       \alpha_2\wedge\alpha_5\wedge\alpha_7+\\
      +&\alpha_3\wedge\alpha_4\wedge\alpha_7-
       \alpha_3\wedge\alpha_5\wedge\alpha_6,
\endsplit
$$
where $\alpha_1,\dots,\alpha_7$ are the basis of the vector space $V$.
{\it If we restrict form of this type to
any $6$-dimensional subspace of $V$ we get a~form of type 2.}

$G_2$ structures are well studied and there is known a~lot of examples
of $G_2$ structures. 

Thus any $G_2$ structure on a~$7$-dimensional manifold
gives a~structure
of type 2 on any $6$-dimensional submanifold. Thus we get a~vast
variety of examples. See for example [J].
\endremark

\end